\documentclass[notitlepage, 11pt]{article}
\usepackage[shortlabels]{enumitem}
\usepackage{amsmath}
\usepackage{amssymb}
\usepackage{amsthm}
\usepackage{abstract}
\usepackage{xcolor}
\usepackage{comment}
\usepackage{mathtools}
\usepackage{graphicx}
\usepackage{caption}
\usepackage{dsfont}
\usepackage{subcaption}
\usepackage{wrapfig}

\usepackage[    backend=bibtex,   style=alphabetic,maxalphanames=4,minalphanames=4 ,maxbibnames=5]{biblatex}
\usepackage{standalone} 
\usepackage{float}
\usepackage{enumitem}
\setlist[itemize]{noitemsep}
\setlist[enumerate]{noitemsep}
\usepackage{hyperref}
\usepackage{cleveref}
\usepackage{tikz}
\usetikzlibrary{arrows.meta}
\usepackage[bottomfloats,belowfloats]{footmisc}

\definecolor{waterfill}{HTML}{DCF2FB}
\definecolor{wateredge}{HTML}{71B7D8}
\definecolor{fishgold}{HTML}{FFBD4C}
\definecolor{fishorange}{HTML}{C47A19}
\definecolor{citecolorhtml}{HTML}{3698c6}

\hypersetup{
    colorlinks=true,
    linkcolor=teal,
    filecolor=black,      
    urlcolor=black,
    citecolor=citecolorhtml,
    }


\newtheorem{theorem}{Theorem}

\usepackage{setspace}
\usepackage[margin=2cm]{geometry} 

\title{\vspace{-2em}On shared responsibility, fish, and the strong law of large numbers}
\date{}
\author{%
  Sahar Diskin\thanks{\textit{Department of Mathematics, ETH Z\"urich, Switzerland. \\Emails: sahardiskinmail@gmail.com, anastasia.kireeva@math.ethz.ch.}}%
  \and Anastasia Kireeva\footnotemark[1]%
}
\begin{document}
\maketitle
\vspace{-2em}
\begin{abstract}
We give a new, short, and intuitive proof of the strong law of large numbers for i.i.d. random variables.
\end{abstract}
While probability theory is largely abstract, many of its fundamental results are easy to observe in real life, even without writing abstract formulas. A prominent example is the \textit{law of large numbers} (LLN). 
While one draw of a random variable, or a result of one experiment may be arbitrarily far from the expected value, it is natural to \emph{expect} that the empirical average of many draws converges to the expected value. The \emph{strong} variant of LLN specifies the mode of this convergence, using one of the strongest notions of convergence in probability theory --- convergence almost surely.

\begin{theorem}\label{thm: main}
Let $X_1,X_2,\ldots$ be i.i.d. real valued random variables with $\mathbb{E}[|X_1|]<\infty$. Let $\mu\coloneqq \mathbb{E}[X_1]$, and for every $n\ge 1$, let $Y_n=(X_1+\cdots+X_n)/n$. Then, $Y_n\to \mu$ almost surely as $n\to \infty$.
\end{theorem}

The statement's intuition can be easily grasped without an extensive mathematical training, and at the same time this law is one of the fundamental results in probability theory. This naturally led to several approaches to prove it: some of the classical treatments include Kolmogorov's truncation argument \cite[Chapter 7]{Feller1971}, backwards martingales \cite[Chapter 4]{Durrett2019}, and the duality of the random walk~\cite{Curien2022}. Here, we give an intuitive proof connecting foundational concepts that are usually already taught in a basic probability course. In fact, we claim that putting rigour aside for a moment, the proof idea is natural for this result.

\begin{figure}[!b]
    \centering
    \begin{subfigure}[c]{0.3\linewidth}
        \centering
        \resizebox{!}{3.8cm}{\begin{tikzpicture}[line cap=round,line join=round]

  \draw[fill=waterfill,draw=wateredge,line width=1pt]
    (0,0) circle[radius=2];
  \node[anchor=north west,inner sep=0pt,font=\Large]
    at (-52:2.15) {$\Omega$};

  \begin{scope}[shift={(0.65,0.65)}]
    \path[fill=fishgold,draw=fishorange,line width=0.8pt]
      (-0.52,0)
      .. controls (-0.52,0.28) and (-0.08,0.34) .. (0.16,0.11)
      -- (0.46,0.30)
      .. controls (0.49,0.18) and (0.39,0.07) .. (0.38,0)
      .. controls (0.39,-0.07) and (0.49,-0.18) .. (0.46,-0.30)
      -- (0.16,-0.11)
      .. controls (-0.08,-0.34) and (-0.52,-0.28) .. (-0.52,0)
      -- cycle;

    \fill[fishorange]
      (-0.35,0.055) circle[radius=0.027];
    \node[font=\large] at (0.53,-0.61) {$A$};
  \end{scope}
\end{tikzpicture}}
        \caption{The orange region $A$ is \\ the event $\{Y_n \not\to \mu\}.\quad $}
        \label{fig:fish-omega}
    \end{subfigure}
    \begin{subfigure}[c]{0.7\linewidth}
        \centering
        \resizebox{!}{4.6cm}{\begingroup%
\definecolor{fishwaterlight}{HTML}{E5F4FC}%
\definecolor{fishwaterblue}{HTML}{C6E5F5}%
\definecolor{fishwateredge}{HTML}{71B7D8}%
\definecolor{fishgoldlight}{HTML}{FFD47F}%
\definecolor{fishgoldorange}{HTML}{FFB343}%
\definecolor{fishoutline}{HTML}{BE781C}%
\definecolor{fishink}{HTML}{31536A}%

\def\FishShape{%
  (-0.52,0)%
  .. controls (-0.52,0.28) and (-0.08,0.34) .. (0.16,0.11)%
  -- (0.46,0.30)%
  .. controls (0.49,0.18) and (0.39,0.07) .. (0.38,0)%
  .. controls (0.39,-0.07) and (0.49,-0.18) .. (0.46,-0.30)%
  -- (0.16,-0.11)%
  .. controls (-0.08,-0.34) and (-0.52,-0.28) .. (-0.52,0)%
  -- cycle%
}%

\def\FishZoomRegion{%
  (-4.18,0.02)%
  .. controls (-4.29,1.77) and (-2.29,2.92) .. (-0.12,2.76)%
  .. controls (2.22,2.94) and (4.17,1.64) .. (4.19,0.07)%
  .. controls (4.33,-1.67) and (2.26,-2.92) .. (0.04,-2.77)%
  .. controls (-2.24,-2.96) and (-4.14,-1.64) .. (-4.18,0.02)%
  -- cycle%
}%

\def\FishBand#1#2#3#4{%
  plot[domain=-3:3,samples=100,variable=\y,smooth]%
    ({#1+0.085*sin(400*\y+#3)},\y)%
  -- plot[domain=3:-3,samples=100,variable=\y,smooth]%
    ({#2+0.085*sin(400*\y+#4)},\y)%
  -- cycle%
}%
\begin{tikzpicture}[line cap=round,line join=round]
  \path[use as bounding box] (-4.45,-3.02) rectangle (4.45,3.02);

  \begin{scope}
    \clip \FishZoomRegion;
    \fill[fishwaterlight] (-4.5,-3.1) rectangle (4.5,3.1);
    \foreach \xl/\xr/\pl/\pr in {
      -2.58/-1.92/10/50,
      -1.26/-0.60/-10/30,
       0.06/0.72/70/0
    }{
      \path[fill=fishwaterblue] \FishBand{\xl}{\xr}{\pl}{\pr};
    }
  \end{scope}
  \draw[draw=fishwateredge,line width=0.8pt,
        dash pattern=on 3.4pt off 2.5pt] \FishZoomRegion;

  \path[fill=fishgoldlight,scale=6] \FishShape;
  \begin{scope}
    \clip[scale=6] \FishShape;
    \foreach \xl/\xr/\pl/\pr in {
      -2.58/-1.92/10/50,
      -1.26/-0.60/-10/30,
       0.06/0.72/70/0
    }{
      \path[fill=fishgoldorange] \FishBand{\xl}{\xr}{\pl}{\pr};
    }
  \end{scope}

  \draw[draw=fishink,line width=0.65pt]
    plot[domain=-1.2835:1.3914,samples=75,variable=\y,smooth]
      ({-2.58+0.085*sin(400*\y+10)},\y);
  \draw[draw=fishink,line width=0.65pt]
    plot[domain=-1.6413:1.6302,samples=75,variable=\y,smooth]
      ({-1.92+0.085*sin(400*\y+50)},\y);
  \draw[draw=fishink,line width=0.65pt]
    plot[domain=-1.7491:1.7461,samples=75,variable=\y,smooth]
      ({-1.26+0.085*sin(400*\y+-10)},\y);
  \draw[draw=fishink,line width=0.65pt]
    plot[domain=-1.6894:1.7070,samples=75,variable=\y,smooth]
      ({-0.60+0.085*sin(400*\y+30)},\y);
  \draw[draw=fishink,line width=0.65pt]
    plot[domain=-1.5254:1.5131,samples=75,variable=\y,smooth]
      ({0.06+0.085*sin(400*\y+70)},\y);
  \draw[draw=fishink,line width=0.65pt]
    plot[domain=-1.0783:1.1165,samples=75,variable=\y,smooth]
      ({0.72+0.085*sin(400*\y+0)},\y);
  \draw[draw=fishoutline,line width=1.05pt,scale=6] \FishShape;

  \node[text=black,font=\large,anchor=base,
        inner sep=0pt,outer sep=0pt] (sliceY) at (-0.73,2.03) {$Y_k>\mu+\varepsilon$};

  \draw[draw=fishink,line width=0.9pt,shorten <=2pt,
        -{Stealth[length=2.3mm,width=1.5mm]}]
    (sliceY.south) -- (-0.93,0.96);
  \node[font=\Large] at (1.95,-1.95) {$A$};
\end{tikzpicture}%
\endgroup%
}
        \caption{The orange slices indexed by $k$ are the events $A\cap \{Y_k>\mu +\varepsilon\}$. If $Y_n \not\to \mu$, there are infinitely many indices $k$ such that $Y_k > \mu + \varepsilon$ for a fixed $\varepsilon > 0$.}
        \label{fig:fish-slices}
    \end{subfigure}
\end{figure}

To say that $Y_n$ converges to $\mu$ almost surely is equivalent to saying that the event that $Y_n$ does not converge to $\mu$ has zero measure. Let us look closely at this region $A$, where $Y_n \not\to \mu$ (see \Cref{fig:fish-omega}). On one hand, an individual experiment $X_1$ has nothing to do with the limit value of $Y_n$ behaving out of place, which mathematically corresponds to $X_1$ being independent of $A$. Hence, on $A$, it still behaves typically and its expected value remains at $\mu$ --- in other words, $X_1$ is \emph{innocent}. However, at the finite \emph{scale}\footnote{We hope the reader will forgive us for this irresistible mischievous wordplay.} things change. If $Y_n$ does not converge to $\mu$, it means that there are infinitely many indices $k$ such that the finite average $Y_k$ is above or below $\mu$. Equivalently, in $A$, there are infinitely many non-empty slices with, say, $Y_k > \mu + \varepsilon$ for a fixed $\varepsilon > 0$,  see \Cref{fig:fish-slices} (these slices may in fact overlap; we will address this slight inconvenience in the formal proof). On each such a slice, since all the experiments are identical and independent, we anticipate the typical value of each one of them --- including $X_1$ --- to be higher than~$\mu$. In other words, they are exchangeable with each other and so \textit{when they all look the same, they all share the blame}. The only way to consolidate this contradiction is if this region $A$ has zero measure, or equivalently, $Y_n$ were to converge almost surely to $\mu$. Below, we give a proof formalising this intuition.

\begin{proof}[Proof of Theorem \ref{thm: main}]
Without loss of generality, let $\mu=0$. For any $\varepsilon> 0$, define events
\[
A_{\varepsilon}\coloneqq \left\{\limsup_{n\to\infty}Y_n>\varepsilon\right\}, B_{\varepsilon} \coloneqq \left\{\liminf_{n\to\infty}Y_n<-\varepsilon\right\}.
\]
Assume there is an $\omega$ such that $Y_n(\omega)\not\to\mu$. Then there exists $\varepsilon>0$ such that for infinitely many indices $n$, either $Y_n(\omega)>\varepsilon$ or $Y_n(\omega) < -\varepsilon$. Therefore, $\omega \in A_\varepsilon \cup B_\varepsilon$, and the statement follows if we show $\mathbb P (A_\varepsilon) = \mathbb P (B_\varepsilon) = 0$. Since the two cases are analogous, we show it for $A_\varepsilon$. Fix such $\varepsilon$, and abbreviate $A\coloneqq A_{\varepsilon}$.

\paragraph{$X_1$ is innocent.} Define $A'=\{\limsup_{n\to\infty} (X_2+\cdots+ X_n)/n>\varepsilon\}$. Since $X_1/n\to 0$ as $n\to\infty$, we have that $\mathds{1}_{A}=\mathds{1}_{A'}$ almost surely. Thus,
\begin{align}
    \mathbb{E}[X_1\mathds{1}_{A}]=\mathbb{E}[X_1\mathds{1}_{A'}]=0\cdot\mathbb{P}(A')=0, \label{eq: innocent}
\end{align}
where the penultimate equality follows since $X_1$ is independent of $A'$.

\paragraph{When they all look the same, $X_1$ shares the blame.} In what remains, we show that $\mathbb{E}[X_1\mathds{1}_{A}]\ge \varepsilon\cdot\mathbb{P}(A)$. Together with \eqref{eq: innocent}, this will imply that $\mathbb{P}(A)=0.$ Note that $A\subseteq \bigcup_{n\ge 1}\{Y_n>\varepsilon\}$, and thus $A\cap \,\bigcup_{k=1}^{N}\{Y_k>\varepsilon\}$ increases to $A$ with $N$. Since $| X_1\mathds{1}_{A\cap \,\bigcup_{k=1}^{N}\{Y_k>\varepsilon\}}|\le| X_1|$ for every $N$, and $\mathbb{E}[|X_1|]<\infty$, by the dominated convergence theorem,
\begin{align}
    \mathbb{E}[X_1\mathds{1}_A]=\mathbb{E}\left[\lim_{N\to\infty}X_1\mathds{1}_{A\cap \,\bigcup_{k=1}^{N}\{Y_k>\varepsilon\}}\right]=\lim_{N\to\infty}\mathbb{E}\left[X_1\mathds{1}_{A\cap \,\bigcup_{k=1}^{N}\{Y_k>\varepsilon\}}\right]. \label{eq: guilty 1}
\end{align}
For every $1\le k \le N$, define $D_{k,N}=\{Y_k>\varepsilon\}\cap\{ Y_{k+1},\ldots, Y_N\le\varepsilon\}.$ Considering these disjoint events allows us to decompose $\bigcup_{k=1}^N\{Y_k>\varepsilon\}$ into $\bigsqcup_{k=1}^ND_{k,N}$. This, in turn, gives
\begin{align}
\mathbb{E}\left[X_1\mathds{1}_{A\cap \,\bigcup_{k=1}^{N}\{Y_k>\varepsilon\}}\right]=\sum_{k=1}^{N}\mathbb{E}\left[X_1\mathds{1}_{A\cap D_{k,N}}\right].\label{eq: guilty 2}    
\end{align}
We will now bound separately each summand on the right hand side. (The motivation to decompose into the disjoint events is that we will bound each summand separately, however even if the expectation is bounded on two events, it is not necessarily that the same bound holds on the union of these events unless they are disjoint\footnote{Consider the expected value of a 6-sided die roll $X$ on the events $E \coloneqq \{\text{the outcome is even}\}$, $F \coloneqq \{\text{the outcome is divisible by }3\}$, and on the union of these two events. We then have that $\mathbb{E}[ X \mathds{1}_{E}] = 4 \ge 4$ and $\mathbb{E} [X \mathds{1}_F] = 4.5 > 4$, and yet $\mathbb{E} [X \mathds 1_{E \cup F}] = 3.5\not \ge4$.}).

Since $X_1,X_2,\ldots$ are i.i.d. and since $\mathds{1}_{A \cap D_{k,N}}$ is symmetric with respect to $X_1,\ldots, X_k$, we have that for every $1\le i \le k$, 
\begin{align*}
\mathbb{E}\left[X_1\mathds{1}_{A\cap D_{k,N}}\right] &= \int x_1 \mathds{1}_{A\cap D_{k,N}} d\mathbb P(x_1, \ldots, x_i,\ldots) = \int x_1 \mathds{1}_{A\cap D_{k,N}} d\mathbb P(x_1) \cdots d\mathbb P(x_i) d\mathbb P(x_{i+1},\ldots) \\
&= \int x_i \mathds{1}_{A\cap D_{k,N}} d\mathbb P(x_i)\cdots d \mathbb P(x_1) d\mathbb P(x_{i+1},\ldots)=\mathbb{E}\left[X_i\mathds{1}_{A\cap D_{k,N}}\right]. 
\end{align*}
Summing over indices $1, \ldots, k$ and dividing by $k$ gives us
\[
\mathbb{E}\left[X_1\mathds{1}_{A\cap D_{k,N}}\right] = \mathbb{E}\left[Y_k\mathds{1}_{A\cap D_{k,N}}\right] \ge \varepsilon\cdot \mathbb{P}(A \cap D_{k,N}).
\]

Substituting the above in \eqref{eq: guilty 2}, we obtain $\mathbb{E}\left[X_1\mathds{1}_{A\cap \,\bigcup_{k=1}^{N}\{Y_k>\varepsilon\}}\right]\ge \varepsilon\cdot \mathbb{P}\left(A\cap \,\bigcup_{k=1}^N\{Y_k>\varepsilon\}\right)$. Therefore, taking $N\to\infty$, by \eqref{eq: guilty 1} we obtain $\mathbb{E}[X_1\mathds{1}_{A}]\ge\varepsilon\cdot \mathbb{P}(A)$, as required.

Considering that $X_1$ is innocent as $\mathbb{E}[X_1 \mathds{1}_A] = 0$, this can only hold if $\mathbb{P}(A) = 0$, as in \Cref{fig:no-fish}.
\par\medskip
\noindent
\begin{minipage}{\linewidth}
    \captionsetup{type=figure}
    \noindent
    \makebox[\linewidth][c]{%
        \raisebox{\depth}{%
            \resizebox{!}{3.8cm}{

\begin{tikzpicture}[line cap=round,line join=round]
  \definecolor{waterfill}{HTML}{DCF2FB}
  \definecolor{wateredge}{HTML}{71B7D8}
  \definecolor{poofink}{HTML}{C47A19}

  \draw[fill=waterfill,draw=wateredge,line width=1pt]
    (0,0) circle[radius=2];
  \node[anchor=north west,inner sep=0pt,font=\Large]
    at (-52:2.15) {$\Omega$};

  \begin{scope}[shift={(0.52,0.55)}]
    \foreach \rayangle/\rayend in {
      0/0.46,45/0.36,90/0.46,135/0.36,
      180/0.46,225/0.36,270/0.46,315/0.36
    }{
      \draw[draw=poofink,line width=0.9pt]
        (\rayangle:0.15) -- (\rayangle:\rayend);
    }
    \node[font=\normalsize,color=poofink] at (0.53,-0.61) {Poof!};
  \end{scope}
\end{tikzpicture}%

}%
        }%
    }%
    \llap{\mbox{\qedhere}}\par
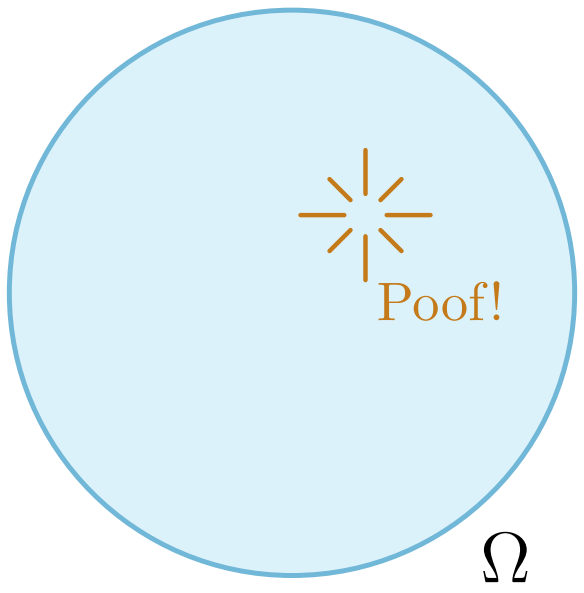
    \caption{No fish (of non-zero measure) in the pond.}
    \label{fig:no-fish}
\end{minipage}
\end{proof}

\paragraph{Acknowledgements} The authors would like to thank Afonso S. Bandeira, Ron Peled, Benny Sudakov, Vincent Tassion, and Maksim Zhukovskii for valuable comments and suggestions. 
\paragraph{Use of AI} AI tools were used to assist with transforming hand-drawn sketches into the tikz code for the figures; no real fish were hurt during the process of human generation of the rest.

\printbibliography
\end{document}